\documentclass[oneside,10pt]{amsart}

\usepackage{amssymb, euscript, longtable}
\newtheorem{Lemma}           {Lemma}
\newtheorem{Corollary}[Lemma]{Corollary}
\newtheorem{Theorem}  [Lemma]{Theorem}

\begin{document}
\title{Projective modules and involutions}
\author{John Murray}
\date{March 11, 2004}
\address{Mathematics Department,
         National University of Ireland - Maynooth,
         Co. Kildare, Ireland.}
\email{jmurray@maths.may.ie}
\begin{abstract}
Let $G$ be a finite group, and let $\Omega:=\{t\in G\mid t^2=1\}$.
Then $\Omega$ is a $G$-set under conjugation.
Let $k$ be an algebraically closed field of characteristic $2$.
It is shown that each projective indecomposable summand of the
$G$-permutation module $k\Omega$ is irreducible and self-dual,
whence it belongs to a real $2$-block of defect zero. This,
together with the fact that each irreducible $kG$-module that
belongs to a real $2$-block of defect zero occurs with multiplicity
$1$ as a direct summand of $k\Omega$, establishes a bijection between
the projective components of $k\Omega$ and the real $2$-blocks of
$G$ of defect zero.
\end{abstract}
\maketitle

Let $G$ be a finite group, with identity element $e$, and let
$\Omega:=\{t\in G\mid t^2=e\}$. Then
$\Omega$ is a $G$-set under conjugation. In this note we describe the
projective components of the permutation module $k\Omega$, where $k$ is
an algebraically closed field of characteristic $2$. By a projective
component we mean an indecomposable direct summand of $k\Omega$ that
is also a direct summand of a free $kG$-module. We show that all such
components are irreducible, self-dual and occur with multiplicity $1$.

This gives an alternative proof of Remark (2) on {p.~254} of \cite{Robinson},
and strengthens Corollaries 3 through 7 of that paper. In addition, we
can give the following quick proof of Proposition 8 in \cite{Robinson}:

\begin{Corollary}
Suppose that $H$ is a strongly embedded subgroup of $G$. Then
$k_H{\uparrow^G}\cong k_G\oplus[\oplus_{i=1}^s P_i]$ where $s\geq 0$
and the $P_i$ are pairwise nonisomorphic self-dual projective irreducible
$kG$-modules.
\begin{proof}
That $H$ is strongly embedded means that $|H|$ is even and $|H\cap H^g|$
is odd, for each $g\in G\backslash H$. Let $t\in H$ be an involution. Then
clearly $C_G(t)\leq H$. So $k_H{\uparrow^G}$ is isomorphic to a submodule of
$(k_{C_G(t)}){\uparrow^G}$. Mackey's theorem implies that every component of
$k_H{\uparrow^G}$, other than $k_G$, is a projective $kG$-module. Being
projective, these modules must be components of $(k_{C_G(t)}){\uparrow^G}$.
The result now follows from Theorem \ref{T:main}.
\end{proof}
\end{Corollary}

Consider the wreath product $G\wr\Sigma$ of $G$ with a cyclic group $\Sigma$
of order $2$. Here $\Sigma$ is generated by an involution $\sigma$ and
$G\wr\Sigma$ is isomorphic to the semidirect product of the base group
$G\times G$ by $\Sigma$. The conjugation action of $\sigma$ on $G\times G$
is given by $(g_1,g_2)^\sigma=(g_2,g_1)$, for all $g_1,g_2\in G$. The elements
of $G\wr\Sigma$ will be written $(g_1,g_2)$, $(g_1,g_2)\,\sigma$ or $\sigma$.

We shall exploit the fact that $kG$ is a $kG\wr\Sigma$-module.
For, as is well-known, $kG$ is an $k{(G\times G)}$-module via:
$x\cdot(g_1,g_2):=g_1^{-1}xg_2$, for each $x\in kG$, and $g_1,g_2\in G$.
The action of $\Sigma$ on $kG$ is induced by the permutation action of $\sigma$
on the distinguished basis $G$ of $kG$: $g^\sigma:=g^{-1}$, for each $g\in G$.
Clearly $\sigma$ acts as an involutary $k$-algebra anti-automorphism
of $kG$. It follows that the actions of $G\times G$ and $\Sigma$ on $kG$ are
compatible with the group relations in $G\wr\Sigma$.

By a \emph{block} of $kG$, or a $2$-block of $G$, we mean an indecomposable
$k$-algebra direct summand of $kG$. Each block has associated to it a
primitive idempotent in $Z(kG)$, a Brauer equivalence class of characters
of irreducible $kG$-modules and a Brauer equivalence class, modulo $2$,
of ordinary irreducible characters of $G$. A block has defect zero if it
is a simple $k$-algebra, and is real if it contains the complex conjugates
of its ordinary irreducible characters. Theorem \ref{T:main} establishes
a bijection between the real $2$-blocks of $G$ that have defect zero and
the projective components of $k\Omega$.

We could equally well work over a complete discrete valaution ring $R$ of
characteristic $0$, whose field of fractions $F$ is algebraically closed,
and whose residue field $R/J(R)$ is $k$. So we use ${\mathcal O}$ to
indicate either of the commutative rings $k$ or $R$.

All our modules are right-modules. We denote the trivial ${\mathcal O}G$-module
by ${\mathcal O}_G$. If $M$ is an ${\mathcal O}G$-module, we use $M{\downarrow_H}$
to denote the restriction of $M$ to $H$. If $H$ is a subgroup of $G$ and $N$ is an
${\mathcal O}H$-module, we use $N{\uparrow^G}$ to denote the induction of $N$ to $G$.
Whenever $g\in G$, we write $\underline{g}$ for $(g,g)\in G\times G$,
and we set $\underline X:=\{\underline x\mid x\in X\}$, for each $X\subset G$.
Other notation and concepts can be found in a standard textbook on modular
representation theory, such as \cite{Alperin} or \cite{NagaoTsushima}. 

If $B$ is a block of ${\mathcal O}G$, then so too is
$B^o=\{x^\sigma\mid x\in B\}$. We call $B$ a real block
if $B=B^o$. Our first result describes the components of
${\mathcal O}G$ as ${\mathcal O}G\wr\Sigma$-module.

\begin{Lemma}\label{L:decomposition}
There is an indecomposable decomposition of
${\mathcal O}G$ as ${\mathcal O}G\wr\Sigma$-module:
$$
{\mathcal O}G=B_1\oplus\ldots\oplus B_r\oplus
(B_{r+1}+B_{r+1}^o)\oplus\ldots\oplus(B_{r+s}+B_{r+s+1}^o).
$$
Here $B_1,\ldots,B_r$ are the real $2$-blocks and
$B_{r+1},B_{r+1}^o,\ldots,B_{r+s},B_{r+s}^o$ are the nonreal $2$-blocks of $G$.
\begin{proof}
This follows from the well-known indecomposable decomposition of
${\mathcal O}G$, as an ${\mathcal O}{(G\times G)}$-module, into a
direct sum of its blocks, and the fact that $B_i^\sigma=B_i$ for
$i=1,\ldots,r$, and $B_{r+j}^\sigma=B_{r+j}^o$ for $j=1,\ldots,s$.
\end{proof}
\end{Lemma}

An obvious but useful fact is that ${\mathcal O}G$ is a permutation module:

\begin{Lemma}\label{L:OGperm}
The ${\mathcal O}G\wr\Sigma$-module ${\mathcal O}G$ is isomorphic
to the permutation module
$({\mathcal O}_{\underline G\times\Sigma}){\uparrow^{G\wr\Sigma}}$.
\begin{proof}
The elements of $G$ form a $G\wr\Sigma$-invariant basis of ${\mathcal O}G$.
Moreover if $g_1,g_2\in G$, then $g_2=g_1\cdot(g_1,g_2)$. So $G$ is a
transitive $G\wr\Sigma$-set. The stabilizer of $e\in{\mathcal O}G$ in
$G\wr\Sigma$ is $\underline G\times\Sigma$. The lemma follows from these facts.
\end{proof}
\end{Lemma}

Let $C$ be a conjugacy class of $G$. Set $C^o:=\{c\in G\mid c^{-1}\in C\}$.
Then $C^o$ is also a conjugacy class of $G$, and $C\cup C^o$ can be regarded
as an orbit of $\underline{G}\times\Sigma$ on the $G\wr\Sigma$-set $G$.
As such, the corresponding permutation module ${\mathcal O}(C\cup C^o)$ is a
${\mathcal O}{\underline G\times\Sigma}$-direct summand of ${\mathcal O}G$.
If $C=C^o$, we call $C$ a real class of $G$. In this case for each $c\in C$
there exists $x\in G$ such that $c^x=c^{-1}$. The point stabilizer of $c$
in $\underline{G}\times\Sigma$ is $\underline{C_G(c)}{<\underline x\sigma>}$.
So ${\mathcal O}C\cong
 ({\mathcal O}_{\underline{C_G(c)}{<\underline x\sigma>}})
 {\uparrow^{\underline G\times\Sigma}}$.
If $C\ne C^o$, we call $C$ a nonreal class of $G$.
In this case the point stabilizer of $c\in C\cup C^o$ in
$\underline{G}\times\Sigma$ is $\underline{C_G(c)}$.
So ${\mathcal O}(C\cup C^o)\cong
({\mathcal O}_{\underline{C_G(c)}}){\uparrow^{\underline G\times\Sigma}}$.

Suppose now that $C_1,\ldots,C_t$ are the real classes of $G$ and that
$C_{t+1},C_{t+1}^o,\ldots,C_{t+u},C_{t+u}^o$ are the nonreal classes.
Then we have:

\begin{Lemma}\label{L:OGpointperm}
There is a decomposition of ${\mathcal O}G$ as an
${\mathcal O}\underline G\times\Sigma$-permutation module:
$$
{\mathcal O}G={\mathcal O}C_1\oplus\ldots\oplus {\mathcal O}C_t\oplus
{\mathcal O}(C_{t+1}\cup C_{t+1}^o)\oplus\ldots\oplus
{\mathcal O}(C_{t+u}\cup C_{t+u+1}^o).
$$
\begin{proof}
This follows from Lemma \ref{L:OGperm} and the discussion above.
\end{proof}
\end{Lemma}

By a quasi-permutation module we mean a direct summand of a permutation
module. Our next result is Lemma 9.7 of \cite{Alperin}. We include a
proof for the convenience of the reader.

\begin{Lemma}\label{L:vertex}
Let $M$ be an indecomposable quasi-permutation ${\mathcal O}G$-module
and suppose that $H$ is a subgroup of $G$ such that $M{\downarrow_H}$
is indecomposable. Then there is a vertex $V$ of $M$ such that $V\cap H$
is a vertex of $M{\downarrow_H}$. If $H$ is a normal subgroup of $G$,
then this is true for all vertices of $M$.
\begin{proof}
Let $U$ be a vertex of $M$. As ${\mathcal O}_U\mid M{\downarrow_U}$ we have
${\mathcal O}_{U\cap H}\mid(M{\downarrow_H}){\downarrow_{U\cap H}}$. But
$U\cap H$ is a vertex of ${\mathcal O}_{U\cap H}$. So Mackey's Theorem implies
that there exists a vertex $W$ of $M{\downarrow_H}$ such that
$U\cap H\leq W$.

As $M{\downarrow_H}$ is a component of the restriction of $M$ to $H$,
Mackey's Theorem shows that there exists $g\in G$ such that $W\leq U^g\cap H$.
Now $U^g$ is a vertex of $M$. So by the previous paragraph, and the uniqueness
of vertices of $M{\downarrow_H}$ up to $H$-conjugacy, there exists $h\in H$
such that $U^g\cap H\leq W^h$. Comparing cardinalities, we see that
$W=U^g\cap H$. So $U^g\cap H$ is a vertex of $M{\downarrow_H}$.

Suppose that $H$ is a normal subgroup of $G$. Then $U\cap H\leq W$
and $W=U^g\cap H=(U\cap H)^g$ imply that $U\cap H=W$.
\end{proof}
\end{Lemma}

R. Brauer showed how to associate to each block of ${\mathcal O}G$
a $G$-conjugacy class of $2$-subgroups, its so-called defect groups.
It is known that a block has defect zero if and only if its defect
groups are all trivial. J. A. Green showed how to associate to each
indecomposable ${\mathcal O}G$-module a $G$-conjugacy class of
$2$-subgroups, its so-called vertices.  He also showed how to
identify the defect groups of a block using its vertices as an
indecomposable ${\mathcal O}(G\times G)$-module.

\begin{Corollary}\label{C:Bvertex}
Let $B$ be a block of ${\mathcal O}G$ and let $D$ be a defect group of $B$.
If $B$ is not real then $\underline D$ is a vertex of $B+B^o$,
as ${\mathcal O}G\wr\Sigma$-module.
If $B$ is real, then there exists $x\in N_G(D)$, with $x^2\in D$,
such that $\underline D<{\underline x\,\sigma}>$ is a vertex of $B$,
as ${\mathcal O}G\wr\Sigma$-module.
In particular, $\Sigma$ is a vertex of $B+B^o$ if and only if
$B$ is a real $2$-block of $G$ that has defect zero.
\begin{proof}
J. A. Green showed in \cite{Green1} that $\underline D$ is a vertex of $B$,
when $B$ is regarded as an indecomposable ${\mathcal O}(G\times G)$-module.
Suppose first that $B$ is not real.
Then $B+B^o=(B{\downarrow_{G\times G}}){\uparrow^{G\wr\Sigma}}$,
for instance by Corollary 8.3 of \cite{Alperin}.
It follows that $B+B^o$ has vertex $\underline D$,
as an indecomposable ${\mathcal O}G\wr\Sigma$-module.

Suppose then that $B=B+B^o$ is real. Lemma \ref{L:OGperm} shows that $B$ is
${\underline G\times\Sigma}$-projective. So we may choose a vertex $V$ of
$B$ such that $V\leq{\underline G}\times\Sigma$. Moreover, $B$ is a
quasi-permutation ${\mathcal O}G\wr\Sigma$-module, and its restriction to
the normal subgroup ${G\times G}$ is indecomposable.
Lemma \ref{L:vertex} then implies that $V\cap({G\times G})=V\cap\underline G$
is a vertex of $B{\downarrow_{{G\times G}}}$. So by Green's result, we may
choose $D$ so that $V\cap\underline G=\underline D$. Now $G\times G$ has
index $2$ in $G\wr\Sigma$. So Green's indecomposability theorem, and the
fact that $B{\downarrow_{{G\times G}}}$ is indecomposable, implies that
$V\not\subseteq{(G\times G)}$. It follows that there exists $x\in N_G(D)$,
with $x^2\in D$, such that $V=\underline D<{\underline x\,\sigma}>$.

If $B$ has defect zero, then $D=<\!e\!>$. So $x^2=e$. In this case,
${<\underline x\,\sigma>}=\Sigma^{(e,x)}$ is $G\wr\Sigma$-conjugate to
$\Sigma$. So $\Sigma$ is a vertex of $B$. Conversely, suppose that $\Sigma$
is a vertex of $B+B^o$. The first paragraph shows that $B$ is a real block
of $G$. Moreover $B$ has defect zero, as $\Sigma\cap\underline G=<e>$.
\end{proof}
\end{Corollary}

We quote the following result of Burry, Carlson and Puig
\cite[4.4.6]{NagaoTsushima} on the Green correspondence:

\begin{Lemma}\label{L:BurryCarlson}
Let $V\leq H\leq G$ be such that $V$ is a $p$-group and $N_G(V)\leq H$.
Let $f$ denote the Green correspondence with respect to $(G,V,H)$.
Suppose that $M$ is an indecomposable ${\mathcal O}G$-module such
that $M{\downarrow_H}$ has a component $N$ with vertex $V$. Then
$V$ is a vertex of $M$ and $N=f(M)$.
\end{Lemma}

We can now prove our main result. Part (ii) is Remark (2) on {p.~254}
of \cite{Robinson}, but our proof is independent of the proof given there.

\begin{Theorem}\label{T:main}
(i) Let $t\in G$, with $t^2=e$. Suppose that $P$ is an indecomposable
projective direct summand of $({\mathcal O}_{C_G(t)}){\uparrow^G}$. Then
$P$ is irreducible and self-dual and occurs with multiplicity $1$ as a
component of $({\mathcal O}_{C_G(t)}){\uparrow^G}$.
In particular $P$ belongs to a real $2$-block of $G$ that has defect zero.
\newline
(ii) Suppose that $M$ is a projective indecomposable ${\mathcal O}G$-module
that belongs to a real $2$-block of $G$ that has defect zero. Then there
exists $s\in G$, with $s^2=e$, such that $M$ is a component of
$({\mathcal O}_{C_G(s)}){\uparrow^G}$. Moreover, $s$ is uniquely determined
up to conjugacy in $G$.
\begin{proof}
If $t=e$ then $P={\mathcal O}_G$. So $P$ is irreducible and self-dual.
The assumption that $P$ is projective and the fact that
${\rm dim}_{\mathcal O}(P)=1$ implies that $|G|$ is odd. So all blocks of
${\mathcal O}G$, in particular the one containing $P$, have defect zero.

Now suppose that $t\ne e$. Let $T$ be the conjugacy class of $G$ that
contains $t$. The permutation module ${\mathcal O}T$ is a direct summand
of the restriction of ${\mathcal O}G$ to $\underline G\times\Sigma$. Regard
$P$ as an ${\mathcal O}\underline G$-module. Let $I(P)$ be the inflation
of this module to $\underline G\times\Sigma$. Then $I(P)$ is a component of
${\mathcal O}T$. As $\Sigma$ is contained in the kernel of $I(P)$, and $P$ is
a projective ${\mathcal O}G$-module, it follows that $I(P)$ has vertex
$\Sigma$ as an indecomposable ${\mathcal O}\underline G\times\Sigma$-module.

By Lemma \ref{L:decomposition}, and the Krull-Schmidt theorem, there exists
a $2$-block $B$ of $G$ such that $I(P)$ is a component of the restriction
$(B+B^o){\downarrow_{\underline G\times\Sigma}}$. An easy computation shows
that $N_{G\wr\Sigma}(\Sigma)=\underline G\times\Sigma$. It then follows from
Lemma \ref{L:BurryCarlson} that $(B+B^o)$ has vertex $\Sigma$ and also that
$I(P)$ is the Green correspondent of $(B+B^o)$ with respect to
$(G\wr\Sigma,\Sigma,\underline G\times\Sigma)$. We conclude from Corollary
\ref{C:Bvertex} that $B$ is a real $2$-block of $G$ that has defect zero.

Let $\hat B$ be the $2$-block of $G\wr\Sigma$ that contains $B$. Then $\hat B$
is real and has defect group $\Sigma$. Let $\hat A$ be the Brauer correspondent
of $\hat B$. Then $\hat A$ is a real $2$-block of $\underline G\times\Sigma$
that has defect group $\Sigma$. Now $\hat A=A\otimes{\mathcal O}\Sigma$, where
$A$ is a real $2$-block of ${\mathcal O}\underline G$ that has defect zero.
In particular $A$ has a unique indecomposable module, and this module is projective,
irreducible and self-dual. Corollary 14.4 of \cite{Alperin} implies that $I(P)$
belongs to $\hat A$. So $P$ belongs to $A$. We conclude that $P$ is irreducible
and self-dual and belongs to a real $2$-block of $G$ that has defect zero.

Now $B$ occurs with multiplicity $1$ as a component of ${\mathcal O}G$,
and $I(P)$ is the Green correspondent of $B$ with respect to
$(G\wr\Sigma,\Sigma,\underline G\times\Sigma)$. So $I(P)$ has multiplicity
$1$ as a component of the restriction of ${\mathcal O}G$ to
$\underline G\times\Sigma$. It follows that $P$ occurs with multiplicity $1$
as a component of $({\mathcal O}_{C_G(t)}){\uparrow^G}$, and with multiplicity
$0$ as a component of $({\mathcal O}_{C_G(r)}){\uparrow^G}$, for $r\in G$ with
$r^2=e$, but $r$ not $G$-conjugate to $t$. This completes the proof of part (i).

Let $R$ be a real $2$-block of $G$ that has defect zero. Then $R$ has vertex
$\Sigma$ as indecomposable ${\mathcal O}G\wr\Sigma$-module. So its Green
correspondent $f(R)$, with respect to $(G\wr\Sigma,\Sigma,\underline G\times\Sigma)$,
is a component of the restriction of ${\mathcal O}G$ to $\underline G\times\Sigma$
that has vertex $\Sigma$.
Lemma \ref{L:OGpointperm} and the Krull-Schmidt theorem imply that
$f(R)$ is isomorphic to a component of ${\mathcal O}(C\cup C^o)$,
for some conjugacy class $C$ of $G$.
Now $\Sigma$ is a central subgroup of $\underline G\times\Sigma$.
So $\Sigma$ must be a subgroup of the point stabilizer of $C\cup C^o$
in $\underline G\times\Sigma$. It follows that $s^2=e$, for each $s\in C$.
Let $N$ denote the restriction of $f(R)$ to $\underline G$, and consider
$N$ as an ${\mathcal O}G$-module. We have just shown that $N$ is a component of
$({\mathcal O}_{C_G(s)}){\uparrow^G}$. Argueing as before, we see that $N$ is an
indecomposable projective ${\mathcal O}G$-module that belongs to a real $2$-block
of $G$ that has defect zero.

The last paragraph establishes an injective map between the real $2$-blocks of $G$
that have defect zero and certain projective components of ${\mathcal O}\Omega$.
As each block of defect zero contains a single irreducible ${\mathcal O}G$-module,
this map must be onto. It follows that the module $M$ in the statement of the
theorem is a component of some permutation module
$({\mathcal O}_{C_G(s)}){\uparrow^G}$, where $s\in G$ and $s^2=e$.
The fact that $s$ is determined up to $G$-conjugacy now follows from the last
statement of the proof of part (i). This completes the proof of part (ii).
\end{proof}
\end{Theorem}

It is possible to simplify the above proof by showing that if $B$ is a real
$2$-block of $G$ that has defect zero, then its Green correspondent, with
respect to $(G\wr\Sigma,\Sigma,\underline G\times\Sigma)$ is $M^{\rm Fr}$,
where $M^{\rm Fr}$ is the Frobenius conjugate of the unique irreducible
${\mathcal O}G$-module that belongs to $B$.

\begin{Corollary}
Let $\Omega=\{t\in G\mid t^2=e\}$. Then there is a bijection between
the real $2$-blocks of $G$ that have defect zero and the projective
components of ${\mathcal O}\Omega$.
\end{Corollary}

Here is a sample application. It was suggested to me by G. R. Robinson.

\begin{Corollary}
Let $n\geq 1$ and let $t$ be an involution in the symmetric group $\Sigma_n$.
If $n=m(m+1)/2$ is a triangular number, and $t$ is a product of
$\lfloor\frac{m^2+1}{4}\rfloor$ commuting transpositions, then there is a
single projective irreducible ${\mathcal O}\Sigma_n$-module, and this module
is the unique projective component of
$({\mathcal O}_{C_{\Sigma_n}(t)}){\uparrow^{\Sigma_n}}$.
For all other values of $n$ or nonconjugate involutions $t$, the modules
$({\mathcal O}_{C_{\Sigma_n}(t)}){\uparrow^{\Sigma_n}}$ are projective free.
\begin{proof}
We give a proof of the following result in \cite[Corollary 8.4]{Murray}:
Let $G$ be a finite group, let $B$ be a real $2$-block of $G$ of defect zero,
and let $\chi$ be the unique irreducible character in $B$. Then there exists a
$2$-regular conjugacy class $C$ of $G$ such that $C=C^o$, $|C_G(c)|$ is odd,
for $c\in C$, and $\chi(c)$ is nonzero, modulo a prime ideal containing $2$.
Moroever, there exists an involution $t\in G$ such that $c^t=c^{-1}$, and for
this $t$ we have $<\chi_{C_G(t)},1_{C_G(t)}>=1$. The existence of $t$
was shown in \cite{Robinson}. The identification of $t$ using the class $C$
was first shown by R. Gow (in unpublished work).

Suppose that $({\mathcal O}_{C_{\Sigma_n}(t)}){\uparrow^{\Sigma_n}}$
has a projective component. Then $\Sigma_n$ has a $2$-block of defect
zero, by Theorem \ref{T:main}. The $2$-blocks of $\Sigma_n$ are indexed
by triangular partitions $\mu=[m,m-1,\ldots,2,1]$, where $m$ ranges over
those natural numbers for which $n-m(m+1)/2$ is even. Moreover, the
$2$-block corresponding to $\mu$ has defect zero if and only if $n=m(m+1)/2$.
In particular, we can assume that $n=m(m+1)/2$, for some $m\geq 1$.

Let $B$ be the unique $2$-block of $\Sigma_n$ that has defect
zero, let $\chi$ be the unique irreducible character in $B$ and
let $g\in\Sigma_n$ have cycle type $\lambda=[2m-1,2m-5,\ldots]$.
Then $|C_{\Sigma_n}(g)|$ is odd. As the parts of $\lambda$ are the
``diagonal hooklengths'' of $\mu$, the Murnaghan-Nakayama formula shows
that $\chi(g)=1$. Now $\lambda$ has $\lfloor(m-1)/2\rfloor$ nonzero parts.
So $g$ is inverted by an involution $t$ that is a product of
$(n-\lfloor(m-1)/2\rfloor)/2=\lfloor\frac{m^2+1}{4}\rfloor$ commuting transpositions.
It follows from Theorem \ref{T:main} and the previous paragraph that the unique
irreducible projective $B$-module occurs with multiplicity $1$ as a component of
$({\mathcal O}_{C_{\Sigma_n}(t)}){\uparrow^{\Sigma_n}}$.
The last statement of the Corollary now follows from Theorem \ref{T:main}.
\end{proof}
\end{Corollary}


\end{document}